.
.
.
\font\sets=msbm10.
\font\stampatello=cmcsc10.
.

\def\spaziolungo{\qquad \qquad \qquad \qquad \qquad \qquad }
\def\1{{\bf 1}}

\def\avesum{\sum_{x\sim N}}

\def\square{\hbox{\vrule\vbox{\hrule\phantom{s}\hrule}\vrule}}
\def\defineq{\buildrel{def}\over{=}}

\def\N{\hbox{\sets N}}

\def\R{\hbox{\sets R}}

\def\modSel{{\widetilde{J}}}

\par
\centerline{\bf A modified Gallagher's Lemma}
\bigskip
\par
\centerline{\stampatello giovanni coppola - maurizio laporta}
\bigskip
{

{\par
{\bf Abstract.} 
First we prove a modified version of the famous Lemma on the mean square estimate for exponential sums, 
by plugging the Cesaro weights 
in the right hand side of Gallagher's inequality. Then we apply it,
in order to establish a mean value estimate for the Dirichlet polynomials.
\footnote{}{\par \noindent {\it Mathematics Subject Classification} : Primary $11{\rm L}07$; Secondary $11{\rm N}37$}
}
\bigskip
\par
\centerline{\bf 1. Introduction and statement of the results.} 
\smallskip
\par
\noindent
Gallagher's Lemma (see [Ga], Lemma 1) is a well-known general mean value estimate for series of the type
$$
S(t)\defineq\sum_{\nu} c(\nu)e(\nu t)\ ,
$$
\par
\noindent
where $e(x)\defineq e^{2\pi ix}$ as usual, 
the frequencies $\nu$ run over a (strictly increasing) sequence of real numbers 
and the coefficients $c(\nu)$ are complex numbers. 
Precisely, it states that, if $S(t)$ is absolutely convergent and $\delta\defineq\theta/T$ with  $\theta\in(0,1)$, then 
$$
\int_{-T}^T|S(t)|^2dt\ll_{\theta}\delta^{-2}\int_{\R}\Big|\sum_{x<\nu \le x+\delta}c(\nu)\Big|^2dx\ .
\leqno{(\star)}
$$
\par
\noindent
Hereafter, according to Vinogradov's notation
$A\ll_{\theta} B$ stands for $|A|\le CB$, where $C>0$ is an unspecified constant that depends on $\theta$, namely $C=C(\theta)$. Typically, in the present context the bounds hold for $T\to 0$ or $T\to \infty$ (more precisely, for $|T|\le T_0$ with a sufficiently small $T_0>0$ or 
 for $T>T_0$ with a sufficiently large $T_0>1$). 
\par
An immediate and renowned consequence is Theorem 1 of [Ga], also known as Gallagher's Lemma for the Dirichlet series.
Indeed, since an absolutely convergent Dirichlet series can be written as
$$
D(t)\defineq\sum_n a_nn^{it}=\sum_{\nu} c(\nu)e(\nu t)\ 
$$
\par
\noindent
by taking $\nu\defineq(2\pi)^{-1}\log n$ and $c(\nu)\defineq a_{e(\nu/i)}$, then, making  the substitution $x=\theta\log y$ 
in $(\star)$
with $\theta\defineq(2\pi)^{-1}$ and recalling that $T=\theta\delta^{-1}$, one immediately has
$$
\int_{-T}^T|D(t)|^2dt\ll 
T^2 \int_0^{+\infty}\Big|\sum_{y<n\le ye^{1/T}}a_n\Big|^2{dy\over y}\ .
\leqno{(\star\star)}
$$
\par
Motivated by our study on the relationship between the Selberg integral
and its modification with the Cesaro weights (see [C] and [CL]), here we give a modified version of the inequality $(\star)$, 
by plugging such weights in
 the right hand side. Indeed, keeping Gallagher's notation, we prove the following variation of his inequality. In passing, we point out that there are further possible generalizations (see the Remark in $\S2$).
\smallskip
\par
\noindent
{\bf Lemma}. {\it Let $T>0$ and $\delta=\theta/T$ with $\theta\in(0,1)$. Then}
$$
\int_{-T}^T|S(t)|^2dt\ll_{\theta}\delta^{-2}\int_{\R}\Big|\sum_{|\nu-x|\le\delta}(1-|\nu-x|\delta^{-1})c(\nu)\Big|^2dx\ . 
\leqno{(\widetilde{\star})}
$$
\par				
\noindent
We remark that, if $S(t)$ is an exponential sum (in  other words, if $c$ has finite support), the integral of the upper bound in $(\star)$ is de facto the {\it Selberg integral} of the arithmetic function $c$ with a vanishing mean value in short intervals (say, $c$ is {\it balanced}), i.e.
$$
J_c(N,\delta)\defineq \avesum \Big| \sum_{x<n\le x+\delta}c(n)\Big|^2\, , 
$$
\par
\noindent
where $x\sim N$ means that $x$ is an integer of the interval $(N,2N]$ with $N\to \infty$. In particular, when $\delta=o(N)$ and $c$ is a real and balanced function, in [CL] (compare $\S4$ up to page 15) it is showed  that one has 
$$
\modSel_c(N,\delta)\defineq \avesum \Big| \sum_{0\le |n-x|\le \delta}\Big( 1-{{|n-x|}\over \delta}\Big)c(n)\Big|^2
\ll J_c(N,\delta)+\delta^3 \Vert c\Vert_{\infty}^2\ , 
$$
\par
\noindent
where $\displaystyle{\Vert c\Vert_{\infty}\defineq \max_{N-\delta<n\le 2N+\delta}|c(n)|}$
and $\modSel_c(N,\delta)$ is the {\it modified} Selberg integral of $c$, that apparently emulates
the integral in the the right hand side of $(\widetilde{\star})$. From this point of view our Lemma can be proposed as a sort of refinement of Gallagher's one. In fact, the above inequality suggests that the modified Selberg integral should be easier to study than  the corresponding Selberg integral. Actually, we expect this to be true, due to the further averaging over the inner short sum  coming from Cesaro weights (compare the discussion in $\S0$ of [CL]  after Theorem 2). 
\par
Noteworthily, with the aid of our Lemma the first author [C] has recently derived a non trivial estimate for the Selberg integral of the three-divisor function $d_3$ under a Conjecture, on the corresponding modified Selberg integral, that is analyzed in [CL].
\par
Here we also obtain a non trivial consequence for Dirichlet polynomials, though it is not as immediate as $(\star\star)$ whenever $(\star)$ is given.
\smallskip
\par
\noindent {\bf Theorem}. {\it For every Dirichlet polynomial $\displaystyle{D(t)=\sum_n a_n n^{it}}$ one has, for} $T\to \infty$, 
$$
\int_{-T}^{T}|D(t)|^2 dt\ll 
T^2 \int_{1}^{+\infty}\Big|\sum_{y-y/T\le n\le y+y/T}\Big(1-{{|n-y|}\over {y/T}}\Big)a_n\Big|^2{dy\over y}
+ \int_{1}^{+\infty}\Big(\sum_{y-\Delta\le n\le y+\Delta}|a_n|\Big)^2{dy\over y} \ ,\leqno(\widetilde{\star\star})
$$
\par
\noindent
{\it where} $\Delta=\Delta(y,T)\defineq y/T+O(y/T^2)$.
\medskip
\par
We think that such a result may be easily generalized to any absolutely convergent series.
\medskip
\par
An application of our Theorem concerns the special case of Dirichlet polynomials approximating Dirichlet series on the critical line $1/2+it$, by taking  $a_n=w(n)b(n)n^{-1/2}$ with a bounded weight $w(n)$ and an
\lq \lq essentially bounded\rq \rq\ arithmetic function
 $b$. Indeed, we prove the following consequence of our Theorem. 
\smallskip
\par
\noindent {\bf Corollary}. {\it Let us consider $P(t)\defineq\displaystyle{\sum_{N_1\le n\le N_2}{{w(n)b(n)}\over {n^{1/2+it}}} }$, where
\thinspace $N_1,N_2$ are positive integers,
$w$ is uniformly bounded and supported in \thinspace $[N_1,N_2]$, and
$|b(n)|\ll_{\varepsilon} n^{\varepsilon}$, $\forall \varepsilon>0$.
Then, as $T\to \infty$, we have
$$
\int_{-T}^{T}|P(t)|^2dt\ll_{\varepsilon} 
T^2 \int_{N_1/2}^{3N_2/2}\Big|\sum_{y-y/T\le n\le y+y/T}\Big(1-{{|n-y|}\over {y/T}}\Big){{w(n)b(n)}\over {n^{1/2}}}\Big|^2{dy\over y}
+ {{N_2^{1+\varepsilon}}\over {T^2}}\ . 
$$
}
\medskip
\par
At least in principle, this Corollary may be applied within the specific situation given in [C0] for the $2k-$th moments of the Riemann $\zeta$ function on the critical line, i.e.$$
I_k(T)\defineq \int_{T}^{2T}\Big|\zeta\Big( {1\over 2}+it\Big)\Big|^{2k} dt\ .
$$
\par				
\noindent
Indeed, Theorem 1.1 of [C0] links $I_k(T)$ to $\modSel_k(N,\delta)$, the Selberg integral of the $k-$divisor function\footnote{$^1$}{Namely, $d_k(n)$ is the number of ways to write $n\ge 1$ as a product of $k$ positive integers and the Dirichlet series generating $d_k(n)$ is the same $\zeta^k(s)$.}
$$
d_k(n)\defineq \sum_{{n_1,\ldots,n_k}\atop {n_1\cdots n_k=n}}1\, ,\quad  (k\in \N)\ .
$$
\par
\noindent
By following [C0] approach, in a forthcoming paper it will be exploited the same link with $\modSel_k(N,\delta)$, where  
essentially it replaces the corresponding Selberg integral in the upper bound for $I_k(T)$, so that our Corollary 
may enter the scene for $b(n)=\widetilde{d_k}(n)$, the balanced part of $d_k(n)$ (compare $\S3$ of [CL]). 
In particular, assuming that the Conjecture given in [CL] holds for the modified Selberg integral of $d_3$, 
such a new approach would lead to the so-called \lq \lq weak sixth moment for the Riemann zeta-function\rq \rq\ (see [CL], $\S8,9$).

\bigskip

\par
\centerline{\bf 2. Proof of the results.} 
\smallskip
\par
\noindent {\stampatello Proof of the Lemma.} 
By following [Ga], where it is introduced the auxiliary function
$$
C_\delta(x)\defineq\sum_{\nu}c(\nu)F_\delta(x-\nu),
$$
\par
\noindent
with $F_\delta(y)=\delta^{-1}$ or $0$ according as $|y|\le\delta/2$ or not, 
analogously we write
$$
\delta^{-2}\int_{\R}\Big|\sum_{|\nu-x|\le\delta}(1-|\nu-x|\delta^{-1})c(\nu)\Big|^2dx=\int_{\R}|\widetilde C_\delta(x)|^2dx\ ,
$$
\par
\noindent
where\enspace $\widetilde C_\delta(x)\defineq\sum_{\nu} c(\nu)\widetilde F_\delta(x-\nu)$\enspace with 
\enspace $\widetilde F_\delta(y)\defineq\max(\delta^{-1}-|y|\delta^{-2},0)$. 
The Fourier transform of $\widetilde C_\delta(x)$ is
$$
\widehat{\widetilde C}_\delta(y)\defineq\int_{\R}\widetilde C_\delta(x)e(-xy)dx=
\sum_{\nu} c(\nu)\int_{\R}\widetilde F_\delta(x-\nu)e(-xy)dx=
$$
$$
=\sum_{\nu} c(\nu)e(-\nu y)\int_{\R}\widetilde F_\delta(t)e(-ty)dt=\widehat{\widetilde F}_\delta(y)S(-y)\ .
$$
\par
\noindent
Thus, by Plancherel's theorem one has
$$
\int_{\R}|\widetilde C_\delta(x)|^2dx=\int_{\R}|\widehat{\widetilde C}_\delta(y)|^2dy=
\int_{\R}|\widehat{\widetilde F}_\delta(y)S(-y)|^2dy
\ ,
$$
\par
\noindent
where the Fourier transform of \enspace $\widetilde F_\delta(x)$ \enspace is given by\footnote{$^2$}{
Note that, as in Gallagher's, we have a normalization condition even in our case, i.e. 
$$
\int_{\R} \widetilde F_\delta(y)dy=\int_{\R} F_\delta(y)dy=1.
$$
}
$$
\widehat{\widetilde F}_\delta(y)=\delta^{-2}\widehat{\max}(\delta-|y|,0)=
\delta^{-2}{{\sin^2(\pi\delta y)}\over{\pi^2y^2}}=\Big({{\sin(\pi\delta y)}\over{\pi\delta y}}\Big)^2=
\widehat{F}_\delta(y)^2\ .
$$
\par
\noindent
Since $\displaystyle{{{\sin(\pi\delta y)}\over{\pi\delta y}}\gg 1}$ when $|y|\le T$, with a constant depending on $\theta$ (recall $\delta=\theta/T$), then 
$$
\int_{-T}^T|S(t)|^2dt=\int_{-T}^T|S(-t)|^2dt 
\ll_{\theta} \int_{\R}|\widehat{\widetilde F}_\delta(y)S(-y)|^2dy
=\delta^{-2}\int_{\R}\Big|\sum_{|\nu-x|\le\delta}(1-|\nu-x|\delta^{-1})c(\nu)\Big|^2dx 
$$
\par				
\noindent
and the Lemma is proved. \hfill $\square$
\medskip
\par
\noindent {\bf Remark.} Evidently, any power of the original transform $\widehat{F}_\delta(y)$ is $\gg 1$ 
(with the implicit constant still depending on $\theta$). However, apart from our case of the square 
that is linked to the Cesaro weights, higher powers would lead to different and much more complicated 
weights. Nevertheless, we are going to explore such implications further in the future. 
\medskip
\par
\noindent {\stampatello Proof of the Theorem.} Let us apply our Lemma to
$$
D(t)=\sum_{\nu} c(\nu)e(\nu t),\enspace \hbox{with}\enspace \nu=(2\pi)^{-1}\log n\ ,\ c(\nu)=a_{e(\nu/i)}\ ,
$$
\par
\noindent
by taking $\theta=(2\pi)^{-1},\ T=\theta\delta^{-1}$ and $x=\theta\log y$. Thus, we get 
$$
\int_{-T}^T|D(t)|^2dt\ll \delta^{-2}\int_{\R}\Big|\sum_{|\nu-x|\le\delta}(1-|\nu-x|\delta^{-1})c(\nu)\Big|^2dx\ll 
$$
$$
\ll 
T^2 \int_0^{+\infty}\Big|\sum_{|\log n-\log y|\le 1/T}(1-T|\log n-\log y|)a_n\Big|^2{dy\over y}\ll 
$$
$$
\ll T^2 \int_{1}^{+\infty}\Big|\sum_{y-(1-1/\tau)y\le n\le y+(\tau-1)y}\Big(1-T\Big|\log\Big(1+{{n-y}\over y}\Big)\Big|\Big)a_n\Big|^2{dy\over y}\ , 
$$
\par
\noindent
where we have set $\tau\defineq e^{1/T}>1$ (see that $\tau\to 1$ as $T\to \infty$, so the $n-$sum is empty for $0<y<1$). 
\smallskip
\par
Since Taylor expansion yields 
$$
y-(1-1/\tau)y=y-{y\over T}+O\Big( {y\over {T^2}}\Big)\ , 
\quad 
y+(\tau-1)y=y+{y\over T}+O\Big( {y\over {T^2}}\Big)\ , 
$$
\par
\noindent
then the Cesaro weight, $\displaystyle{1-T|\log(1+{{n-y}\over y})|}$, is bounded for the present range of $n$, while in both ranges 
$$
0\le |n-(y-y/T)|\ll y/T^2
\enspace \hbox{\rm and} \enspace 
0\le |n-(y+y/T)|\ll y/T^2
$$
\par
\noindent
we have 
$$
1-T\Big|\log\Big(1+{{n-y}\over y}\Big)\Big|\ll {1\over T}\ .
$$
\par
\noindent
Accordingly we write 
$$
\int_{-T}^T|D(t)|^2dt\ll 
T^2 \int_{1}^{+\infty}\Big|\sum_{y-y/T\le n\le y+y/T}\Big(1-T\Big|\log\Big(1+{{n-y}\over y}\Big)\Big|\Big)a_n\Big|^2{dy\over y} + 
$$
$$
+ \int_{1}^{+\infty}\Big(\sum_{0\le |n-(y-y/T)|\ll y/T^2}|a_n|+\sum_{0\le |n-(y+y/T)|\ll y/T^2}|a_n|\Big)^2{dy\over y} \ .
$$
\par
\noindent
Then $(\widetilde{\star\star})$ follows since by Taylor expansion again we have 
$$
T\Big|\log\left(1+{{n-y}\over y}\right)\Big| - {{|n-y|}\over {y/T}}\ll 
{{T(n-y)^2}\over {y^2}}\ll 
{1\over T}
$$
\par
\noindent
for \enspace $y-y/T\le n\le y+y/T$ \enspace and this yields
$$
T^2\int_{1}^{+\infty}\Big|\sum_{y-y/T\le n\le y+y/T}\Big(1-T\Big|\log\Big(1+{{n-y}\over y}\Big)\Big|\Big)a_n\Big|^2{dy\over y}\ll 
$$
$$				
\ll T^2 \int_{1}^{+\infty}\Big|\sum_{y-y/T\le n\le y+y/T}\Big(1-{{|n-y|}\over {y/T}}\Big)a_n\Big|^2{dy\over y} + 
 \int_{1}^{+\infty}\Big(\sum_{-\Delta\le n-y\le \Delta}|a_n|\Big)^2{dy\over y} \ ,
$$
\par
\noindent
where recall that $\Delta=y/T+O(y/T^2)$\ .\hfill \square
\bigskip
\par
\noindent {\stampatello Proof of the Corollary.} 
Let us apply the Theorem to $D(-t)$ with $a_n=w(n)b(n)n^{-1/2}$
and write
$$
\int_{-T}^{T}|P(t)|^2 dt\ll_{\varepsilon} 
T^2 \int_{1}^{+\infty}\Big|\sum_{y-y/T\le n\le y+y/T}\Big(1-{{|n-y|}\over {y/T}}\Big){{w(n)b(n)}\over {n^{1/2}}}\Big|^2{dy\over y}
+ N_2^{\varepsilon}\int_{1}^{+\infty}\Big(\sum_{y-\Delta\le n\le y+\Delta}{1\over {\sqrt n}}\Big)^2{dy\over y}\ . 
$$
\par
\noindent
Since $w$ has support  in $[N_1,N_2]$, then
$$
\int_{-T}^{T}|P(t)|^2 dt\ll_{\varepsilon} 
T^2 \int_{N_1/2}^{3N_2/2}\Big|\sum_{y-y/T\le n\le y+y/T}\Big(1-{{|n-y|}\over {y/T}}\Big){{w(n)b(n)}\over {n^{1/2}}}\Big|^2{dy\over y}
+ N_2^{\varepsilon}\int_{N_1/2}^{3N_2/2}{{\Delta(y,T)^2}\over {y^2}} dy\ll_{\varepsilon} 
$$
$$
\ll_{\varepsilon} T^2 \int_{N_1/2}^{3N_2/2}\Big|\sum_{y-y/T\le n\le y+y/T}\Big(1-{{|n-y|}\over {y/T}}\Big){{w(n)b(n)}\over {n^{1/2}}}\Big|^2{dy\over y}
+ N_2^{1+\varepsilon}{1\over {T^2}}\ , 
$$
\par
\noindent
using \enspace $\Delta=\Delta(y,T)\ll y/T$,\enspace whence the Corollary. \hfill $\square$ 

\bigskip
\par
\noindent {\bf Acknowledgement}. The authors wish to thank Alberto Perelli for helpful comments. 

\bigskip

\par
\centerline{\stampatello References}
\medskip
\item{\bf [C0]} Coppola, G. \thinspace - \thinspace {\sl On the Selberg integral of the $k$-divisor function and the $2k$-th moment of the Riemann zeta-function} \thinspace - \thinspace Publ. Inst. Math. (Beograd) (N.S.) {\bf 88(102)} (2010), 99--110. $\underline{\tt MR\enspace 2011m\!:\!11173}$  \thinspace - \thinspace available online 
\item{\bf [C]} Coppola, G. \thinspace - \thinspace {\sl On the Selberg integral of the three-divisor function $d_3$} \thinspace - \thinspace available online at the address http://arxiv.org/abs/1207.0902 (see version 3)
\item{\bf [CL]} Coppola, G. and Laporta, M. \thinspace - \thinspace {\sl Generations of correlation averages} \thinspace - \thinspace http://arxiv.org/abs/1205.1706 (see version 3)
\item{\bf [Ga]} Gallagher, P. X. \thinspace - \thinspace {\sl A large sieve density estimate near $\sigma =1$} \thinspace - \thinspace Invent. Math. {\bf 11} (1970), 329--339. $\underline{\tt MR\enspace 43\# 4775}$ 

\bigskip

\par
\leftline{\tt Giovanni Coppola\spaziolungo \spaziolungo \qquad \qquad \enspace \thinspace Maurizio Laporta}
\leftline{\tt Universit\`a degli Studi di Salerno\spaziolungo \thinspace Universit\`a degli Studi di Napoli}
\leftline{\tt Home address \negthinspace : \negthinspace Via Partenio \negthinspace 12 \negthinspace -\spaziolungo Dipartimento di Matematica e Appl.}
\leftline{\tt - 83100, Avellino(AV), ITALY\spaziolungo \qquad \qquad \qquad \qquad \enspace \thinspace Compl.Monte S.Angelo}
\leftline{\tt e-page : $\! \! \! \! \! \!$ www.giovannicoppola.name\qquad \qquad \qquad \qquad \qquad \quad \thinspace Via Cinthia - 80126, Napoli, ITALY}
\leftline{\tt e-mail : gcoppola@diima.unisa.it\spaziolungo \qquad \enspace \thinspace e-mail : mlaporta@unina.it}

\bye